# On totally umbilical isometric immersions with a special kind of Liouville's metric in the $\mathbb{R}^3$ space

*Alexander Pigazzini*

ABSTRACT. *The aim of this short article is to investigate the possibility of existence of totally umbilical isometric immersions in $\mathbb{R}^3$ with isothermal parametrization and harmonic metric.*

## INTRODUCTION

In this paper we investigate a class of surfaces in $\mathbb{R}^3$. Precisely, we deal with sufficiently smooth surfaces which have, in isothermal coordinates, an harmonic metric coefficient and an "isothermal" second fundamental form, which means that its matrix is pointwise a multiple of the identity. In particular, those kind of surfaces are totally umbilical ones, hence, as it is well known we are reduced to planes or spheres respectively if the gaussian curvature is zero or not.
Actually our approach is more general and it does not take into account this fact. Surfaces of our interest are a special kind of so called Liouville's surfaces which can be endowed with a so called Liouville's metric, i.e. that admit a nontrivial quadratic first integral of its geodesic flow. Such metrics can be indeed put in the local form $ds^2 = h(x,y)(dx^2 + dy^2)$ where $h > 0$ satisfies $h_{xx} + h_{yy} = 0$ (see **[1]** for details).
We prove that the only surface as required is the plane and the proof can be obtained exploiting the classical system of Gauss and Codazzi-Mainardi (see **[2]** as reference) coupled with the assumption of harmonicity of the metric. As a corollary we say that the sphere cannot be parametrized by this special kind of Liouville's metric.

**Proposition.** *A totally umbilical surface, in isothermal coordinates, has necessarily $l_{ij}(u,v) = e(u,v)\delta_{ij}$ as coefficients of the second fundamental form.*

*Proof. T*he Weingarten's matrix, in isothermal coordinates, is

(1) $$W(v) = \begin{vmatrix} -\dfrac{l_{11}}{g_{11}} & -\dfrac{l_{12}}{g_{11}} \\ -\dfrac{l_{12}}{g_{11}} & -\dfrac{l_{22}}{g_{11}} \end{vmatrix}$$

and the eigenvalues are the principal curvatures. Then, since S is totally umbilical, it must be $\lambda_1 = \lambda_2$, where $\lambda_1$ and $\lambda_2$ are solutions of

(2) $$\left(\lambda + \frac{l_{11}}{g_{11}}\right)\left(\lambda + \frac{l_{22}}{g_{11}}\right) - \frac{l_{11}^2}{g_{11}^2} = g_{11}^2 \lambda^2 + g_{11}\lambda l_{22} + l_{11}l_{22} - l_{12}^2 = 0.$$

Therefore

(3) $$l_{11}^2 + l_{22}^2 - 2l_{11}l_{22} + 4l_{12}^2 = 0$$

Namely $(l_{11} - l_{22})^2 = -4l_{12}^2$ which has only solutions of type $l_{ij}(u,v) = e(u,v)\delta_{ij}$.

∎

**Lemma.** *Let $U$ be open and connected in $\mathbb{R}^2$ and let $i: U \to \mathbb{R}^3$ be an immersion. Suppose that $g_{ij}(u,v) = E(u,v)\delta_{ij}$, where $E(u,v)$ is a strictly positive and harmonic function, and $l_{ij}(u,v) = e(u,v)\delta_{ij}$ as coefficients of the second fundamental form of $i$. Then $e(u,v) = 0$.*

*Proof. Let us write the* system of Gauss and Codazzi-Mainardi equations. For convenience we will use $g_{ij}$ and $l_{ij}$ instead of $g_{ij}(u,v)$ and $l_{ij}(u,v)$. Remember that in isothermal coordinates we have the following equalities of Christoffel symbols:

(4) $$\Gamma_{11}^1 = \Gamma_{12}^2 = -\Gamma_{22}^1; \quad \Gamma_{22}^2 = \Gamma_{12}^1 = -\Gamma_{11}^2.$$

Considering $g_{ij}(u,v) = E(u,v)\delta_{ij}$ and the harmonicity of $E(u,v)$, the system of Gauss and Codazzi-Mainardi equations is reduced to:

(5) $$\begin{cases} \dfrac{\partial l_{22}}{\partial u} - \dfrac{\partial l_{12}}{\partial v} = \left[\dfrac{1}{2}\dfrac{\partial g_{11}}{\partial u} g_{11}^{-1}\right](l_{11} + l_{22}) \\ \dfrac{\partial l_{11}}{\partial v} - \dfrac{\partial l_{12}}{\partial u} = \left[\dfrac{1}{2}\dfrac{\partial g_{11}}{\partial v} g_{11}^{-1}\right](l_{11} + l_{22}) \\ \dfrac{1}{2} g_{11}^{-1}|\nabla g_{11}|^2 = l_{11}l_{22} - l_{12}^2. \end{cases}$$

Since $l_{ij}(u,v) = e(u,v)\delta_{ij}$ the system (5) becomes

(6) $\begin{cases} \frac{\partial l_{11}}{\partial u} = \left[\frac{\partial g_{11}}{\partial u} g_{11}^{-1}\right] l_{11} \\ \frac{\partial l_{11}}{\partial v} = \left[\frac{\partial g_{11}}{\partial v} g_{11}^{-1}\right] l_{11} \\ (l_{11})^2 (g_{11}^{-1})^2 = \frac{1}{2} |\nabla g_{11}|^2 (g_{11}^{-1})^3. \end{cases}$

From the first two of (6) we obtain

(7) $l_{11} = c g_{11}$, with $c$ an arbitrary constant.

Substituting (7) in the third of (6) we obtain:

(8) $2c^2 g_{11}^3 = |\nabla g_{11}|^2$.

We now prove that from (8) we have only two cases:

- $c = 0; g_{11} = constant$
- $c \neq 0; g_{11} = 0$.

First, if $c = 0$, then clearly $g_{11}$ must be constant because from (8) we have $\frac{\partial g_{11}}{\partial u} = \frac{\partial g_{11}}{\partial v} = 0$.

Assume now that $c \neq 0$.

Differentiating (8) with respect to $u$ and $v$ we obtain:

$$\begin{cases} \frac{\partial g_{11}}{\partial u}\frac{\partial^2 g_{11}}{\partial u^2} + \frac{\partial g_{11}}{\partial v}\frac{\partial^2 g_{11}}{\partial u \partial v} = 3c^2 g_{11}^2 \frac{\partial g_{11}}{\partial u} \\ \frac{\partial g_{11}}{\partial u}\frac{\partial^2 g_{11}}{\partial u \partial v} + \frac{\partial g_{11}}{\partial v}\frac{\partial^2 g_{11}}{\partial v^2} = 3c^2 g_{11}^2 \frac{\partial g_{11}}{\partial v} \end{cases}$$

and, coupled with the harmonic equation $\frac{\partial^2 g_{11}}{\partial u^2} + \frac{\partial^2 g_{11}}{\partial v^2} = 0$, we deduce that

$$\begin{cases} \frac{\partial g_{11}}{\partial u}\frac{\partial^2 g_{11}}{\partial u^2} + \frac{\partial g_{11}}{\partial v}\frac{\partial^2 g_{11}}{\partial u \partial v} = 3c^2 g_{11}^2 \frac{\partial g_{11}}{\partial u} \\ \frac{\partial g_{11}}{\partial u}\frac{\partial^2 g_{11}}{\partial u \partial v} + \frac{\partial g_{11}}{\partial v}\frac{\partial^2 g_{11}}{\partial v^2} = 3c^2 g_{11}^2 \frac{\partial g_{11}}{\partial v} \\ \frac{\partial^2 g_{11}}{\partial u^2} + \frac{\partial^2 g_{11}}{\partial v^2} = 0. \end{cases}$$

Hence

$$
(9) \begin{cases} \dfrac{\partial^2 g_{11}}{\partial u^2} = 3c^2 g_{11}{}^2 \dfrac{(\frac{\partial g_{11}}{\partial u})^2 - (\frac{\partial g_{11}}{\partial v})^2}{(\frac{\partial g_{11}}{\partial u})^2 + (\frac{\partial g_{11}}{\partial v})^2} \\[2mm] \dfrac{\partial^2 g_{11}}{\partial u \partial v} = \dfrac{6 c^2 g_{11}{}^2 \frac{\partial g_{11}}{\partial u}\frac{\partial g_{11}}{\partial v}}{(\frac{\partial g_{11}}{\partial u})^2 + (\frac{\partial g_{11}}{\partial v})^2} \\[2mm] \dfrac{\partial^2 g_{11}}{\partial v^2} = 3c^2 g_{11}{}^2 \dfrac{(\frac{\partial g_{11}}{\partial v})^2 - (\frac{\partial g_{11}}{\partial u})^2}{(\frac{\partial g_{11}}{\partial u})^2 + (\frac{\partial g_{11}}{\partial v})^2}. \end{cases}
$$

Using (9) to expand the identity $\dfrac{\partial^3 g_{11}}{\partial u^2 \partial v} - \dfrac{\partial^3 g_{11}}{\partial u \partial v \partial u} = 0$, we find that

$$
(10) \qquad g_{11} \frac{\partial g_{11}}{\partial v} [3c^2 g_{11}{}^3 - \left(\frac{\partial g_{11}}{\partial u}\right)^2 - \left(\frac{\partial g_{11}}{\partial v}\right)^2] = 0
$$

while the identity $\dfrac{\partial^3 g_{11}}{\partial u \partial v^2} - \dfrac{\partial^3 g_{11}}{\partial v^2 \partial u} = 0$ yields

$$
(11) \qquad g_{11} \frac{\partial g_{11}}{\partial u} [3c^2 g_{11}{}^3 - \left(\frac{\partial g_{11}}{\partial u}\right)^2 - \left(\frac{\partial g_{11}}{\partial v}\right)^2] = 0.
$$

Suppose, by contradiction, that $g_{11}$ is not identically zero. Let $(u_0, v_0) \in U$ be such that $g_{11}(u_0, v_0) \neq 0$. Since (8) we also can say that $\dfrac{\partial g_{11}}{\partial u}(u_0, v_0)$ and $\dfrac{\partial g_{11}}{\partial v}(u_0, v_0)$ cannot be simultaneously zero. Therefore (10) and (11) give $3c^2 g_{11}(u_0, v_0)^3 = |\nabla g_{11}(u_0, v_0)|^2$ which contradicts (8).

Since now by assumption $g_{11} > 0$, we are reduced to the case $c = 0$, and then (7) yields the conclusion.

∎

**Remark 1.** In reference **[3]**, the Ricci tensor in isothermal coordinates is given by $-2R_{ij} = \Delta g_{ij} + Q(g^{-1}, \partial g)$, where $Q(g^{-1}, \partial g)$ denotes a sum of terms which are quadratic in the metric inverse $g^{-1}$ and its first derivatives $\partial g$. From this we can see that:
If $\Delta g_{ij} < 0$ then $R_{ij} > 0$
If $\Delta g_{ij} = 0$ then $R_{ij} \geq 0$
but if $\Delta g_{ij} = 0$ with $l_{ij}(u,v) = e(u,v)\delta_{ij}$ then, by our *Lemma*, we have only $R_{ij} = 0$.

∎

**Theorem.** *The only immersion $i: U \to \mathbb{R}^3$, with $g_{ij}(u,v) = E(u,v)\delta_{ij}$, where $E(u,v)$ is a strictly positive and harmonic function, and with $l_{ij}(u,v) = e(u,v)\delta_{ij}$ as coefficients of the second fundamental form, is a plane.*

*Proof.* Let $i$ be the immersion of a plane, then choosing the orthogonal Cartesian coordinates, we have $g_{ij}(u,v) = \delta_{ij}$ as coefficients of the first fundamental form and $l_{ij}(u,v) = 0$ as coefficients of the second fundamental form. Viceversa, if some immersion satisfies $g_{ij}(u,v) = E(u,v)\delta_{ij}$, where $E(u,v)$ is a strictly positive and harmonic function, and $l_{ij}(u,v) = e(u,v)\delta_{ij}$ as coefficients of the second fundamental form, as seen in *Lemma,* we have that $e(u,v) = 0$. Then $i$ is affine; in other words $i(U)$ is (contained in) a plane.

∎

**Remark 2.** It is well know that the only totally umbilical surfaces in $\mathbb{R}^3$ are the planes and the spheres, respectively if the gaussian curvature is zero or not. Taking into account this fact, the proof of the *Theorem* could be simplified; we leave open this possible kind of simplification.

**Corollary.** Considering the *Proposition*, that tells us that every isometric and totally umbilical immersion, in isothermal coordinates, has necessarily second fundamental form as $l_{ij}(u,v) = e(u,v)\delta_{ij}$, and considering the *Theorem,* we can say, that does not exist a possible parametrization for the spheres that is compatible with this special kind of Liouville's metric.

*E-mail address:* pigazzinialexander18@gmail.com